\documentclass[12pt]{article}
\usepackage{amsfonts}
\usepackage{latexsym}
\newtheorem{theorem}{Theorem}

\begin{document}

\title{Random sampling in chirp space}
\author{Eric Carlen\thanks{%
School of Mathematics, Georgia Institute of Technology, Atlanta, GA 30332,
USA, carlen@math.gatech.edu, Work partially 
supported by U.S. National Science Foundation 
 grant DMS 03-00349}, R. Vilela Mendes\thanks{%
Centro de Matem\'{a}tica e Aplica\c{c}\~{o}es Fundamentais, Av. Gama Pinto
2, 1649-003 Lisboa, Portugal, vilela@cii.fc.ul.pt;
http://label2.ist.utl.pt/vilela/ }\thanks{%
Centro de Fus\~{a}o Nuclear - EURATOM/IST Association, Instituto Superior
T\'{e}cnico, Av. Rovisco Pais 1, 1049-001 Lisboa, Portugal} }
\date{}
\maketitle

\begin{abstract}
For the space of functions that can be approximated by linear chirps, we
prove a reconstruction theorem by random sampling at arbitrary rates.
\end{abstract}


\section{Introduction. Chirps, fractional Fourier transform and random
sampling}

In the Fr\'{e}chet space of infinitely differentiable functions, with a
topology defined by a family of norms 
\begin{equation}
\left\| f\left( t\right) \right\| _{n,k}=\sup_{t\in \Bbb{R}}\left|
t^{n}f^{\left( k\right) }\left( t\right) \right| \qquad n,k=0,1,2,\cdots
\label{1.1}
\end{equation}
the Fourier operator $\mathcal{F}_{1}$%
\begin{equation}
\left( \mathcal{F}_{1}f\right) \left( \omega \right) =F\left( \omega \right)
=\frac{1}{\left( 2\pi \right) ^{1/2}}\int_{-\infty }^{\infty }e^{-i\omega
t}f\left( t\right) dt  \label{1.2}
\end{equation}
is an invertible transform with 
\begin{equation}
f\left( t\right) =\frac{1}{\left( 2\pi \right) ^{1/2}}\int_{-\infty
}^{\infty }e^{i\omega t}F\left( \omega \right) d\omega  \label{1.3}
\end{equation}

In the time-frequency plane with orthogonal axis $t$ and $\omega $, the
Fourier transform has a geometrical interpretation as a rotation of the
signal by $\alpha _{1}=\frac{\pi }{2}$. The fractional Fourier transform 
\cite{Ozaktas} \cite{Bultheel1} \cite{Bultheel2} corresponds also to a
rotation in $t-\omega $ plane, but now by a non-integer multiple of $\frac{%
\pi }{2}$, $\alpha _{b}=b\frac{\pi }{2}$, namely 
\begin{equation}
\left( \mathcal{F}_{b}f\right) \left( \zeta \right) =F_{b}\left( \zeta
\right) =\frac{e^{-\frac{i}{2}\left( \textnormal{sgn}\left( \sin \alpha
_{b}\right) \frac{\pi }{2}-\alpha _{b}\right) }}{\left( 2\pi \left| \sin
\alpha _{b}\right| \right) ^{1/2}}\int_{-\infty }^{\infty }e^{\left( -i\frac{%
t\zeta }{\sin \alpha _{b}}+\frac{i}{2}\cot \alpha _{b}\left( t^{2}+\zeta
^{2}\right) \right) }f\left( t\right) dt  \label{1.4}
\end{equation}
with inverse $\mathcal{F}_{-b}$.

Another way to look at the Fourier transform equation (\ref{1.2}) is as a
decomposition of $f$ into a combination of harmonics (the basis functions $%
\left\{ e^{-i\omega t}\right\} $). As a result, the Fourier transform is a
convenient way to code the signal when this one is a superposition of a
(small) number of harmonics.

Instead of $e^{-i\omega t}$ the kernel of the fractional Fourier transform (%
\ref{1.4}) is a linear chirp $e^{-i\left( \omega -ct\right) t}$ with $\omega
=\frac{\zeta }{\sin \alpha _{b}}$ and $c=\frac{1}{2}\cot \alpha _{b}$. This
suggests that a more general basis to expand a signal, with arbitrary
features in the $\omega -t$ plane, is a basis of linear chirps.

The choice of the appropriate basis is an important issue when
reconstructing signals from limited data and in particular for the
reconstruction from non-uniform sampling. Sampling at irregular intervals is
sometimes imposed by the nature of the phenomenon being observed as in the
case of geological measurement, anemometry or radar signals. Reconstruction
of the signal from such measurements poses special problems, because the
usual spectral estimation methods are designed to work with equidistantly
spaced time series. Special methods have therefore been developed to deal
with the problem of irregular sampling of band-limited functions\cite{Landau}
\cite{Seip} \cite{Benedetto} \cite{Feichtinger} \cite{Grochenig} \cite{Say}.

Another problem, in the reconstruction or interpolation of signals, in the
band-limited space, occurs when the sampling rate is below the Nyquist rate,
especially in view of the Beurling-Landau density theorem\cite{Landau} \cite
{Grochening2}. In this case irregular sampling of an appropriate type,
instead of being a nuisance, may be of help for the asymptotically exact
reconstruction of signals. This occurs in the space $\mathcal{A}$ of almost
periodic functions and the basic mathematical result is the following:

\begin{theorem}
(Collet\cite{Collet}) Let $x_{n}=n\lambda +X_{n}$ with $X_{n}$ being a
sequence of i.i.d. random variables uniformly distributed in [$0$,$\lambda $%
]. Then, almost every configuration \{$x_{n}$\} of the point process has the
property that if $f$ is any complex almost periodic function satisfying 
\begin{equation}
f\left( x_{n}\right) =0\hspace{3cm}\forall n\in \Bbb{Z}  \label{1.5}
\end{equation}
Then $f\equiv 0$.
\end{theorem}

\begin{center}
------------
\end{center}

Recall that $f$ is almost periodic\cite{Besicovitch} \cite{Corduneanu} if it
is uniformly continuous and if $\forall \varepsilon >0$ $\exists \Lambda
\left( \varepsilon \right) $ such that any interval 
\begin{equation}
\lbrack a,a+\Lambda \left( \varepsilon \right) ]  \label{1.6}
\end{equation}
contains a number $\tau \left( a,\varepsilon \right) $ such that 
\begin{equation}
\sup_{x}\left| f\left( x+\tau \left( a,\varepsilon \right) \right) -f\left(
x\right) \right| \leq \varepsilon  \label{1.7}
\end{equation}

Given $f\left( x\right) $ almost periodic, $\forall \varepsilon >0$ $\exists 
$ a trigonometric polynomial $g$ approximating $f$ uniformly, that is $%
\exists B_{1}\left( \varepsilon \right) ,...,B_{a}\left( \varepsilon \right) 
$ such that 
\begin{equation}
g\left( x\right) =\sum_{k=1}^{a}B_{k}e^{i2\pi \omega _{k}x}  \label{1.8}
\end{equation}
\begin{equation}
\sup_{x\in \Bbb{R}}\left| f\left( x\right) -g\left( x\right) \right| \leq
\varepsilon  \label{1.9}
\end{equation}

The approximation by trigonometric polynomials and Collet's theorem provides
a basis for asymptotically exact reconstruction algorithms at rates much
below Nyquist's rate. However if the local signal frequencies
vary in time the approximation by trigonometric polynomials is not very
convenient. A more general \textit{chirp basis} would be more appropriate.
Here one provides a generalization of Collet's theorem to a space of
functions that may be approximated by linear chirps. It contains, as a
particular case, the space of almost periodic functions.

\section{A uniform approximation result for random sampling in chirp space}

Instead of (\ref{1.7}), the approximation by trigonometric polynomials (Eqs. 
\ref{1.8} and \ref{1.9}) provides an alternative, and equivalent,
characterization of almost periodic functions. Likewise we define the space $%
\mathcal{LC}$ of \textit{linear chirp functions} as the space of functions $%
f $ such that $\forall \varepsilon >0$ $\exists $ a finite number of real
number sets $\left( \omega _{1},c_{1},\alpha _{1},B_{1}\right) ,...,\left(
\omega _{k},c_{k},\alpha _{k},B_{k}\right) $ such that 
\begin{equation}
g\left( x\right) =\sum_{j=1}^{k}B_{j}e^{i\left\{ \omega _{j}+c_{j}\left(
x-\alpha _{j}\right) \right\} x}  \label{2.1}
\end{equation}
and 
\begin{equation}
\sup_{x\in \Bbb{R}}\left| f\left( x\right) -g\left( x\right) \right| \leq
\varepsilon  \label{2.2}
\end{equation}
(the space of almost periodic functions corresponds to the $c_{j}=0$ case).

The space $\mathcal{LC}$ of linear chirp functions is strictly larger than
the space of almost periodic functions. It suffices to consider $e^{ix^{2}}$%
. If it were an almost periodic functions, it would exist $\xi $ such that 
\[
\left| e^{i\left( x+\xi \right) ^{2}}-e^{ix^{2}}\right| <\varepsilon \hspace{%
2cm}-\infty <x<\infty 
\]
with $\varepsilon <1$. Choose $x^{*}=\frac{k\pi -\xi ^{2}}{2\xi }$ and $k$
odd. Then 
\[
\left| e^{i\left( x^{*}+\xi \right) ^{2}}-e^{ix^{*2}}\right| =2>\varepsilon 
\]
a contradiction.

For functions in $\mathcal{LC}$ one obtains the following:

\begin{theorem}
\textit{Let }$x_{n}=n\lambda +X_{n}$\textit{\ with }$X_{n}$\textit{\ being a
sequence of i.i.d. random variables uniformly distributed in [}$0$\textit{,}$%
\lambda $\textit{]. Then, almost every configuration \{}$x_{n}$\textit{\} of
the point process has the property that if }$f$\textit{\ is a function in
the linear chirp space satisfying }
\[
f\left( x_{n}\right) =0\hspace{3cm}\forall n\in \Bbb{Z},
\]
\textit{then }$f\equiv 0$\textit{.}
\end{theorem}

For the proof one needs the following :

\begin{theorem}
\textbf{\ }\textit{For almost every configuration }$\left\{ x_{n}\right\} $%
\textit{\ of the random process, one has} 
\[
\lim_{L\rightarrow \infty }\frac{1}{2L+1}\sum_{-L\leq n\leq L}e^{i\left(
\omega x_{n}+cx_{n}^{2}\right) }=0
\]
\textit{for real }$\omega $\textit{\ and }$c$\textit{\ with }$\omega \neq 0$%
\textit{.}
\end{theorem}

The proof may follow similar steps as the proof of proposition 5 in \cite
{Collet}. A simpler argument uses the invariant measure properties of the
random dynamical system in the circle 
\begin{equation}
y_{n}=\frac{\omega }{2\pi }\left( n\lambda +X_{n}\right) +\frac{c}{2\pi }%
\left( n\lambda +X_{n}\right) ^{2}\hspace{2cm}\left( \textnormal{mod}1\right)
\label{L1}
\end{equation}

Let first $c=0,\omega \neq 0$ and (without loss of generality) $\lambda =1$.
What the dynamical system (\ref{L1}) does is to cover the circle with
intervals of length $\frac{\omega }{2\pi }$ and, in each one, to choose a
point at random. Because $X_{n}$ has uniform distribution in each interval,
the distribution of $y_{n}$ is also uniform. Therefore, by the ergodic
theorem 
\begin{equation}
\lim_{L\rightarrow \infty }\frac{1}{2L+1}\sum_{-L\leq n\leq L}e^{i\omega
x_{n}}=\left\langle e^{i2\pi y_{n}}\right\rangle _{S^{1}}=0  \label{L2}
\end{equation}
for generic sequences \{$X_{n}$\}.

For $c$ and $\omega \neq 0$, the proof is more delicate. At each step $n$,
of the random dynamical system \ref{L1}, the interval where the random point
is chosen is 
\begin{equation}
\left[ \frac{\omega }{2\pi }n+\frac{c}{2\pi }n^{2},\frac{\omega }{2\pi }%
\left( n+1\right) +\frac{c}{2\pi }\left( n+1\right) ^{2}\right]   \label{L3}
\end{equation}
For all 
\[
n>n^{*}=\frac{2\pi -\omega }{2c}-\frac{1}{2}
\]
the interval in (\ref{L3}) wraps one or more times over the circle. In the
interval (\ref{L3}) the distribution of the random points would be
\[
\rho _{n}\left( y\right) =\frac{2\pi }{\sqrt{\omega ^{2}+8\pi cy}}
\]
When the interval wraps around the circle the maximum possible deviation
from uniformity of the density is
\begin{equation}
\Delta \rho _{n}=\frac{2\pi }{\sqrt{\omega ^{2}+4c\left( \omega
n+cn^{2}\right) }}  \label{L4}
\end{equation}

Therefore we have a sequence $\{y_{n}\}$ of independent random variables
with values in $[0,1]$ and each one has a density $\rho _{n}(y)$ with
\[
\rho _{n}(y)=1+r_{n}(y)
\]
where 
\[
\sup_{x\in [0,1]}|r_{n}(x)|=a_{n}\ .
\]

We are interested in the case in which 
\begin{equation}
a_{n}= {\cal O}(1/n)  \label{L5}
\end{equation}
and what to check whether 
\begin{equation}
\lim_{N\to \infty }{\frac{1}{N}}\sum_{j=1}^{N}\sin (2\pi jX_{j})=0\ ,
\label{L6}
\end{equation}
and also to relate the speed of convergence in the second limit to the speed
of convergence in the first limit.

For each $n$, define $Y_{n}$ by 
\begin{equation}
Y_{n}=\sin (2\pi X_{n})  \label{L7}
\end{equation}
and let $\mu _{n}$ be its mean. Define $S_{N}$ by 
\begin{equation}
S_{N}=\sum_{n=1}^{N}Y_{n}\ .  \label{L8}
\end{equation}

Also, define $M_{n}(\xi )$ by 
\begin{equation}
M_{n}(\xi )=E\left( e^{\xi Y_{n}}\right) \ .  \label{L9}
\end{equation}
By Jensen's inequality, 
\[
M_{n}(\xi )\ge e^{\xi \mu _{n}}\ .
\]

The following lemma is a simple variant of a well known estimate of Cramer.
The variation, probably not new, is that our sequence of random variables, $%
\{Y_{n}\}$, is not identically distributed. \medskip 

\noindent \textbf{Lemma 1} \textit{Suppose that for each $N$ there is a
convex function $\Phi _{N}(\xi )$ such that 
\begin{equation}
\sum_{n=1}^{N}\ln M_{n}(\xi /N)\le \Phi _{N}(\xi )\ .  \label{L10}
\end{equation}
Let $I_{N}(y)$ denote the Legendre transform of $\Phi _{N}$; i.e., 
\begin{equation}
I_{N}(y)=\sup_{\xi }\left( y\xi -\Phi _{N}(\xi )\right) \ .  \label{L11}
\end{equation}
Also, let 
\begin{equation}
m_{N}={\frac{1}{N}}\sum_{n=1}^{N}\mu _{n}\ .  \label{L12}
\end{equation}
}

\textit{Then, for any number $y>m_{N}$, 
\begin{equation}
\ln \left( \mathrm{Pr}\left( S_{N}/N>y\right) \right) \le -I_{N}(y)\ .
\label{L13}
\end{equation}
} \medskip 

\noindent \textbf{Proof:} By the Chebychev
inequality, for all $\xi \ge 0$, 
\[
\mathrm{Pr}\left( S_{N}/N>y\right) \le e^{-y\xi }E\left( e^{\xi
S_{N}/N}\right) \ .
\]
But by the independence, 
\[
E\left( e^{\xi S_{N}/N}\right) =\prod_{n=1}^{N}M_{n}(\xi /N)\ ,
\]
and so 
\begin{eqnarray}
{\ln \left( \mathrm{Pr}\left( S_{N}/N>y\right) \right) } &\le &{-\xi
y+\sum_{n=1}^{N}\ln M_{n}(\xi /N)}  \label{L15} \\
&&  \nonumber
\end{eqnarray}
Since $\ln M_{n}(\xi )\ge \xi \mu _{n}$ for each $n$ and $\xi $, it follows $%
m_{N}\xi \le \Phi _{N}(\xi )$. Then
\[
y\xi -\Phi _{N}(\xi )\le y\xi -m_{N}\xi \le 0
\]
for each $\xi <0$ and $y>m_{N}$. Then 
\[
I_{N}(y)=\sup_{\xi \ge 0}\left( y\xi -\Phi _{N}(\xi )\right) \ .
\]
Therefore, (\ref{L13}) follows from (\ref{L15}).\medskip 

The same sort of reasoning shows that for any number $y<m_{N}$, 
\[
\ln \left( \mathrm{Pr}\left( S_{N}/N<y\right) \right) \le -I_{N}(y)\ 
\]

We now apply this in the following setting. From (\ref{L5}) we may deduce
estimates on each $\ln M_{n}(\xi )$ of the form 
\[
\ln M_{n}(\xi )\le \mu _{n}\xi +C\xi ^{2}\ .
\]
Then, 
\[
\sum_{n=1}^{N}\ln M_{n}(\xi )\le N(m_{N}\xi +C\xi ^{2})\ .
\]

In this case, we can apply the Lemma with 
\[
\Phi _{N}(\xi )=m_{N}\xi +{\frac{C}{N}}\xi ^{2}\ .
\]

Computing the Legendre transform, we then have 
\[
I_N(y) = {\frac{N}{4C}}(m_N - y)^2\ .
\]

Now, define $b_{N}$ by 
\[
b_{N}^{2}={\frac{8C\ln N}{N}}\ .
\]
Then 
\[
\mathrm{Pr}\left( S_{N}/N>m_{N}+b_{N}\right) \le e^{-Nb_{N}^{2}/4C}=N^{-2}\ .
\]
Since this is summable, the Borel--Cantelli lemma implies that the probability for $S_{N}/N$ to exceed $m_{N}+b_{N}$
infinitely often is zero. This large deviation estimate implies (\ref{L6}).
A similar reasoning applies for the cosine series. This completes the proof
of theorem 3.\medskip 

\noindent\textbf{Proof of theorem 2}\textit{:}

Let $f\left( x_{n}\right) =0$, $\forall n\in \Bbb{Z}$ and $g\left( x\right) $
be its $\varepsilon -$approximation by linear chirp polynomials. Then 
\begin{eqnarray*}
&&\left| \lim_{L\rightarrow \infty }\frac{1}{2L+1}\sum_{-L\leq n\leq
L}e^{-i\left\{ \omega x_{n}+cx_{n}^{2}\right\} }g\left( x_{n}\right) \right|
\\
&=&\left| \lim_{L\rightarrow \infty }\frac{1}{2L+1}\sum_{-L\leq n\leq
L}e^{-i\left\{ \omega x_{n}+cx_{n}^{2}\right\} }\left( g\left( x_{n}\right)
-f\left( x_{n}\right) \right) \right| \leq \varepsilon
\end{eqnarray*}
for all $\omega $ and $c$. Inserting now Eq.(\ref{2.1}) in the left-hand
side of the above equation one obtains 
\[
\left| \lim_{L\rightarrow \infty }\frac{1}{2L+1}\sum_{-L\leq n\leq
L}\sum_{j=1}^{k}B_{j}e^{i\left\{ \left( \omega _{j}-c_{j}\alpha _{j}-\omega
\right) x_{n}+\left( c_{j}-c\right) x_{n}^{2}\right\} }\right| \leq
\varepsilon 
\]
Choosing $\omega =\omega _{j}-c_{j}\alpha _{j}$ , $c=c_{j}$ and using the
lemma, one concludes that for almost every configuration $\left\{
x_{n}\right\} $, 
\[
\left| B_{j}\right| \leq \varepsilon 
\]
for all $j$ in the linear chirp approximation.

Because this result holds for all $\varepsilon $ and the linear chirp basis
functions are kernels to the fractional Fourier transform, one concludes
that the function $f$ has zero fractional Fourier spectrum. Therefore it is
the zero function.$\Box $

As in the case of functions in the almost periodic space, the
above result may be used to estimate functions in the linear chirp space by
random sampling. If from a time series $h\left( x_{n}\right) $, one obtains,
by the appropriate algorithm, a linear chirp approximation $g\left( x\right) 
$ coinciding with the sampled function on a typical sequence $\left\{
x_{n}\right\} $, that is 
\[
f\left( x_{n}\right) =g\left( x_{n}\right) -h\left( x_{n}\right) =0 
\]
then, in the above defined space, one knows that $g\left( x\right) =h\left(
x\right) $ for all $x$.

\section{Nonlinear chirps} 

By a nonlinear chirp, we mean a linear combination of functions of the form
$e^{ip(x)}$ where $p$ is a real polynomial of some degree $m>2$. As before, let $c$ denote the leading coefficient.
The analysis leading to Theorem 3 easily extends to include nonlinear as well as linear chirps.  Although these
nonlinear chirp functions have no connection with the fractional Fourier transform when, it is of interest to
see how far one may extend the analysis, especially as the extension requires very little modification.

To see how this goes, observe that for large $|x|$, the leading term in $p(x)$ dominates, and so there
is an $N< \infty$ so that for all $|x|>N$,  
\begin{equation}\label{nlc}
|p'(x)| \ge (|c|m/2)|x|^{m-1} > 2\pi\ .
\end{equation}
 In particular, for 
$|n|>N$, $p(x)$ is monotone on $[n,n+1]$.

Let $\sigma_n(y)$ be the denisty for the random variable $y = p(x)$ with $x$
chosen uniformly in $[n,n+1]$ for $|n|>N$.  We may assume without loss of generality
that $p(x)$ is increasing on this interval, and in this case,
$$\sigma_n(y) = {1\over p'(x(y))}\ ,$$
where $p(n) \le y \le p(n+1)$, and $x(y)$ is the unique solution of $p(x) = y$. 

As before, the interval $[p(n),p(n+1)]$ ``wraps'' one or more times around the circle, inducing the density
$\rho_n(y)$ on the circle. The variation between the maximum and the minimum of $\rho_n$
is no more than the maximum value of $\sigma_n(y)$. By our assumption (\ref{nlc}), this is no more
 than $2m|c|/n^{m-1}$.
 
 Therefore, we have that
 $$\rho_n(y) = 1 + r_n(y)$$
 with 
 $$\sup_{x\in [0,1]}|r_{n}(x)|=a_{n} = {\cal O}(n^{m-1})\ .$$
 This is an improvement over (\ref{L5}), and from here, the proof of the analog of Theorem 3
 proceeds as before.
 
 \medskip
 
\noindent{\bf Acknowledgement: } E.C. would like to thank Centro de Matem\'{a}tica e Aplica\c{c}\~{o}es Fundamentais for hospitality
 during his work in Lisbon on this paper.

\end{document}